\documentclass[11pt]{article}
\usepackage[utf8]{inputenc}
\usepackage{graphicx}
\usepackage{amsmath, amsthm, amssymb} 
\usepackage{hyperref} 
\usepackage{multicol}
\usepackage{multirow}
\usepackage{caption}
\usepackage{xcolor}
\usepackage{subcaption}
\usepackage[labelfont=bf]{caption}
\usepackage{pdflscape}
\usepackage[noend]{algpseudocode}

\usepackage{algorithm}

\newcommand{\myalgsheader}[0]

\algnewcommand{\IIf}[1]{\State\algorithmicif\ #1\ \algorithmicthen}
\algnewcommand{\EndIIf}{\unskip\ \algorithmicend\ \algorithmicif}
\algnewcommand{\IElse}[1]{\State\algorithmicelse\ #1\ \algorithmicthen}
\algnewcommand{\IfThenElse}[3]{
  \State \algorithmicif\ #1\ \algorithmicthen\ #2\ \algorithmicelse\ #3}

\begin{document}
\title{Mixed Precision GMRES-based Iterative Refinement with Recycling\footnote{We acknowledge funding from the Charles University PRIMUS project
No. PRIMUS/19/SCI/11, Charles University Research Program No. UNCE/SCI/023, and the Exascale Computing Project (17-SC-20-SC), a collaborative effort of the U.S. Department of Energy Office of Science and the National Nuclear Security Administration.}}

\author{Eda Oktay and Erin Carson}
\date{}
\maketitle

\paragraph{Abstract.}
With the emergence of mixed precision capabilities in hardware, iterative refinement schemes for solving linear systems $Ax=b$ have recently been revisited and reanalyzed in the context of three or more precisions. These new analyses show that under certain constraints on condition number, the LU factorization of the matrix can be computed in low precision without affecting the final accuracy. Another promising technique is GMRES-based iterative refinement, which, in contrast to the standard approach, use GMRES preconditioned by the low-precision triangular factors to solve for the approximate solution update in each refinement step. This more accurate solution method extends the range of problems which can be solved with a given combination of precisions. However, in certain settings, GMRES may require too many iterations per refinement step, making it potentially more expensive than simply recomputing the LU factors in a higher precision. 

Krylov subspace recycling is a well-known technique for reusing information across sequential invocations of a Krylov subspace method on systems with the same or a slowly changing coefficient matrix. In this work, we incorporate the idea of Krylov subspace recycling into a mixed precision GMRES-based iterative refinement solver. The insight is that in each refinement step, we call preconditioned GMRES on a linear system with the same coefficient matrix $A$, with only the right-hand side changing. In this way, the GMRES solves in subsequent refinement steps can be accelerated by recycling information obtained from the first step. We perform extensive numerical experiments on various random dense problems, Toeplitz problems (prolate matrices), and problems from real applications, which confirm the benefits of the recycling approach.

\section{Introduction} \label{sec:first_intro}
Iterative refinement (IR) a commonly-used approach for solving nonsingular linear systems $Ax=b$, where $A \in \mathbb{R}^{n \times n}, x,b \in \mathbb{R}^n$ \cite{w:63}. In each refinement step $i$, the approximate solution ${x}_{i+1} \in \mathbb{R}^n$ is updated by adding a correction term obtained from the previously computed approximate solution ${x}_{i}$. The initial approximate solution ${x}_0$ is usually computed by Gaussian elimination with partial pivoting (GEPP). The computed ${L}$ and ${U}$ factors of $A$ are then reused to solve the correction equation $A {d}_{i+1} = {r}_i$, where ${r}_i=b-A{x}_i$ is the residual. Finally, the original approximate solution is refined by the correction term, ${x}_{i+1} = {x}_{i} + {d}_{i+1}$, and the process continues iteratively until the desired accuracy is achieved.

Recently, mixed precision hardware has become commonplace in supercomputing architecture. This has inspired the development of a new mixed precision benchmark for supercomputers, called HPL-AI, which is based on mixed precision iterative refinement and on which today's top supercomputers exceed exascale performance; see, e.g., \cite{hplai,k:20,top500}. Algorithm \ref{alg:sir} shows a general IR scheme in a mixed precision setting. The authors in \cite{ch:18} used three hardware precisions in the algorithm: $u_f$ is used for LU factorization, $u_r$ for residual calculation, and $u$ (the working precision) for the remaining computations. The authors also introduce a fourth quantity, $u_s$, which is the effectively precision of the solve (taking on the values of $u$ or $u_r$ depending on the particular solver used in line \ref{corrsolve}). We assume throughout this work that $u_f\geq u\geq u_r$. 

\begin{algorithm}[htbp!]
	\caption{General Iterative Refinement Scheme \label{alg:sir}} 
	\begin{algorithmic}[1]
		\Require{$n \times n$ matrix $A$; right-hand side $b$; maximum number of refinement steps $i_{max}$.}
		\Ensure{Approximate solution $x_{i+1}$ to $Ax = b$.}
		\State{Compute LU factorization $A \approx LU$ in precision $u_f$}.
		\State{Solve $Ax_0 = b$ by substitution in precision $u_f$; store $x_0$ in precision $u$.}
		\For{$i = 0$ : $i_{max} - 1$}    
		\State{Compute $r_i = b - Ax_i$ in precision $u_r$; store in precision $u$.}
		\State{Solve $Ad_{i+1} = r_i$ in precision $u_s$; store $d_{i+1}$ in precision $u$. \label{corrsolve}}
		\State{Compute $x_{i+1} = x_i + d_{i+1}$ in precision $u$.}
		\IIf{converged} return $x_{i+1}$. \EndIIf
		\EndFor
	\end{algorithmic}
\end{algorithm}

For a given combination of precisions and choice of solver, it is well-understood under which conditions Algorithm \ref{alg:sir} will converge and what the limiting accuracy will be. The constraint for convergence is usually stated via a constraint on the infinity-norm condition number of the matrix $A$; see, e.g., \cite{ch:18}. In the case that the computed LU factors are used to solve for the correction in line \ref{corrsolve}, often referred to as ``standard IR'' (SIR), then $\kappa_\infty(A)= \Vert A \Vert_\infty \Vert A^{-1} \Vert_\infty$ must be less than $u_f^{-1}$ in order for convergence to be guaranteed. 

To relax this constraint on condition number, the authors of \cite{ch:17} and \cite{ch:18} devised a mixed precision GMRES-based iterative refinement scheme (GMRES-IR), given in Algorithm \ref{alg:gmresir}. In GMRES-IR, the correction equation is solved via left-preconditioned GMRES, where the computed LU factors of $A$ are used as preconditioners. This results in a looser constraint on condition number in order to guarantee the convergence of forward and backward errors; in the case that the preconditioned matrix is applied to a vector in each iteration of GMRES in double the working precision, we require $\kappa_\infty(A)\leq u^{-1/2}u_f^{-1}$, and in the case that a uniform precision is used within GMRES, we require $\kappa_\infty(A)\leq u^{-1/3}u_f^{-2/3}$; see \cite{h:21}. If these constraints are met, we are guaranteed that the preconditioned GMRES method will converge to a backward stable solution after $n$ iterations and that the iterative refinement scheme will converge to its limiting accuracy. We note that in order to guarantee backward stability, all existing analyses (e.g., \cite{ch:17, ch:18, h:21}) assume that unrestarted GMRES is used within GMRES-IR.

\begin{algorithm}[htbp!]
	\caption{GMRES-IR \cite{ch:17} \label{alg:gmresir}}
	\begin{algorithmic}[1]
		\Require{$n \times n$ matrix $A$; right-hand side $b$; maximum number of refinement steps $i_{max}$; GMRES convergence tolerance $\tau$.}
		\Ensure{Approximate solution ${x}_{i+1}$ to $Ax = b$.}
		\State{Compute LU factorization $A \approx LU$ in precision $u_f$.}
		\State{Solve $Ax_0 = b$ by substitution in precision $u_f$; store $x_0$ in precision $u$.}
		\For{$i = 0$ : $i_{max} - 1$}   
		\State{Compute $r_i = b - Ax_i$ in precision $u_r$; scale $r_i = r_i/\|r_i\|_\infty$; store in precision $u$.}
		\State{Solve $U^{-1}L^{-1}Ad_{i+1}=U^{-1}L^{-1}r_i$ by GMRES with tolerance $\tau$ in working precision $u$, with matrix-vector products with $\tilde{A}=U^{-1}L^{-1}A$ computed in precision $u^2$; store $d_{i+1}$ in precision $u$. \label{corrsolveg}}
		\State{Compute $x_{i+1}=x_i+\|r_i\|_\infty d_{i+1}$ in precision $u$.}
		\IIf{converged} return $x_{i+1}$. \EndIIf
		\EndFor
	\end{algorithmic}
\end{algorithm}



We also note that existing analyses are unable to say anything about how fast GMRES will convergence in each refinement step, only that it will do so within $n$ iterations. 
However, if indeed $n$ iterations are required to converge in each GMRES solve, this can make GMRES-IR more expensive than simply computing the LU factorization in higher precision and using SIR. Unfortunately, GMRES convergence is incredibly difficult to predict. In fact, for any set of prescribed eigenvalues, one can construct a linear system for which GMRES will stagnate entirely until the $n$th iteration \cite{gps:96}. The situation is better understood at least in the case of normal matrices, see, e.g., \cite{lt:04}. The worst-case scenario in the case of normal $A$ is when eigenvalues are clustered near the origin, which can cause complete stagnation of GMRES \cite{lt:04}. After the preconditioning step in GMRES-IR, all eigenvalues of the preconditioned matrix $(U^{-1}L^{-1}A)$ ideally become 1 in the absence of finite precision error in computing LU and within GMRES. However, in practice, since we have inexact LU factors, if $A$ has a cluster of eigenvalues near the origin, this imperfect preconditioner may fail to shift some of them away from the origin, which can cause GMRES to stagnate. For instance, when random dense matrices having geometrically distributed singular values are used in the multistage iterative refinement algorithm devised in \cite{oc:21}, the authors showed that for relatively large condition numbers relative to precision $u_f$, GMRES tends to perform $n$ iterations in each refinement step. 

\begin{figure}[h!]
	\centering
	\includegraphics[trim={3cm 8cm 4cm 8cm},clip, width=.45\textwidth]{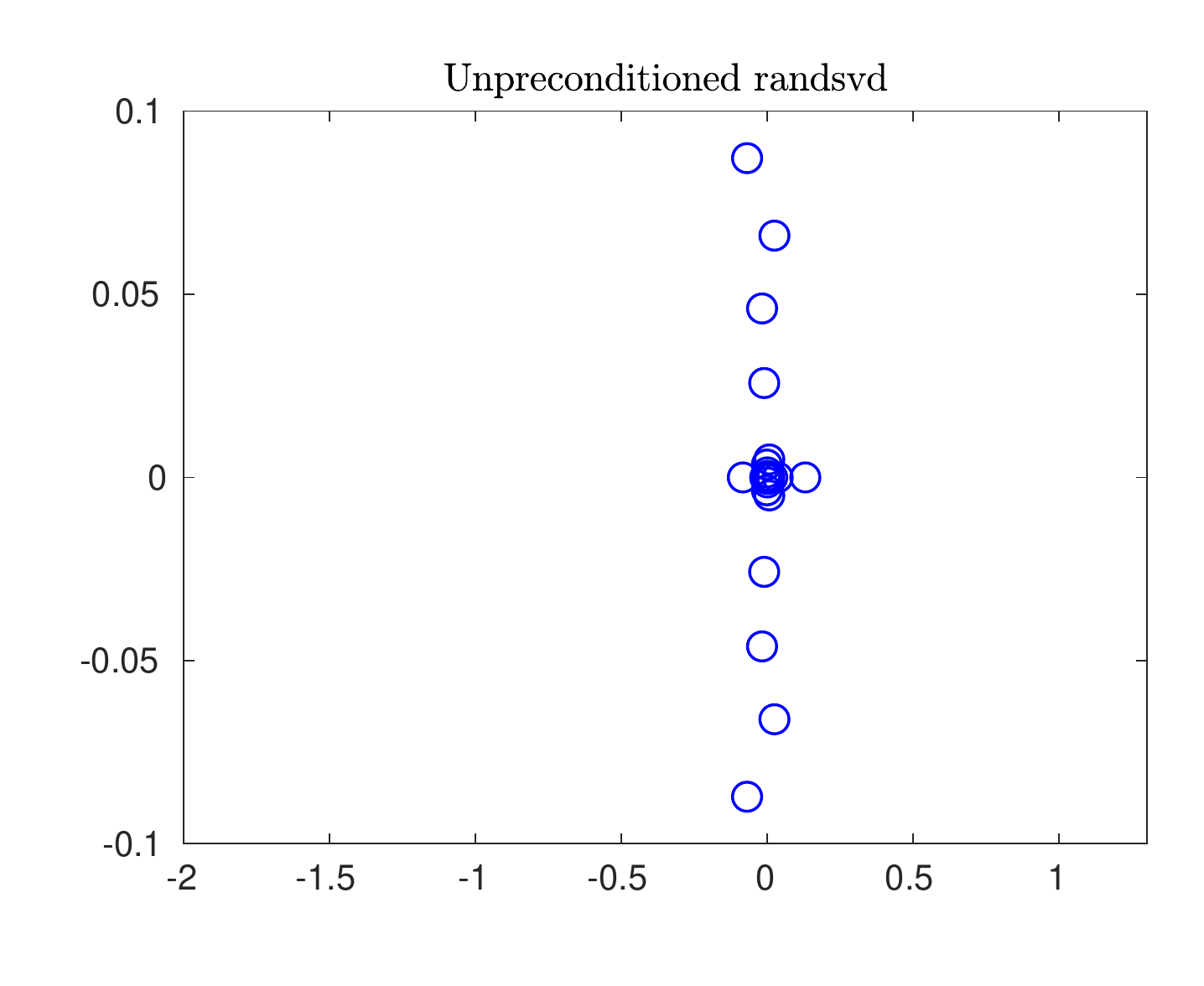}
	\includegraphics[trim={3cm 8cm 4cm 8cm},clip, width=.45\textwidth]{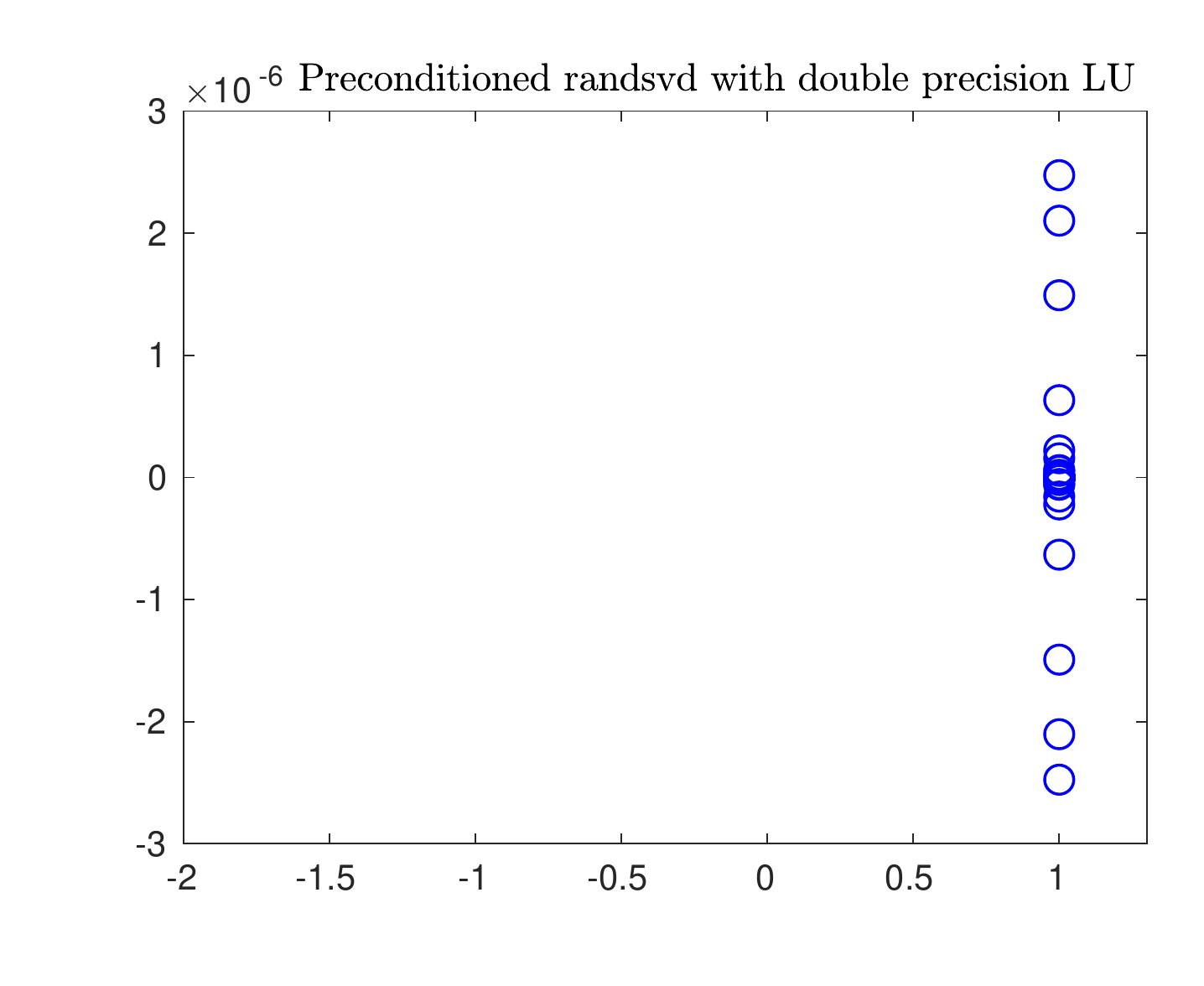}
	\includegraphics[trim={3cm 8cm 4cm 8cm},clip, width=.45\textwidth]{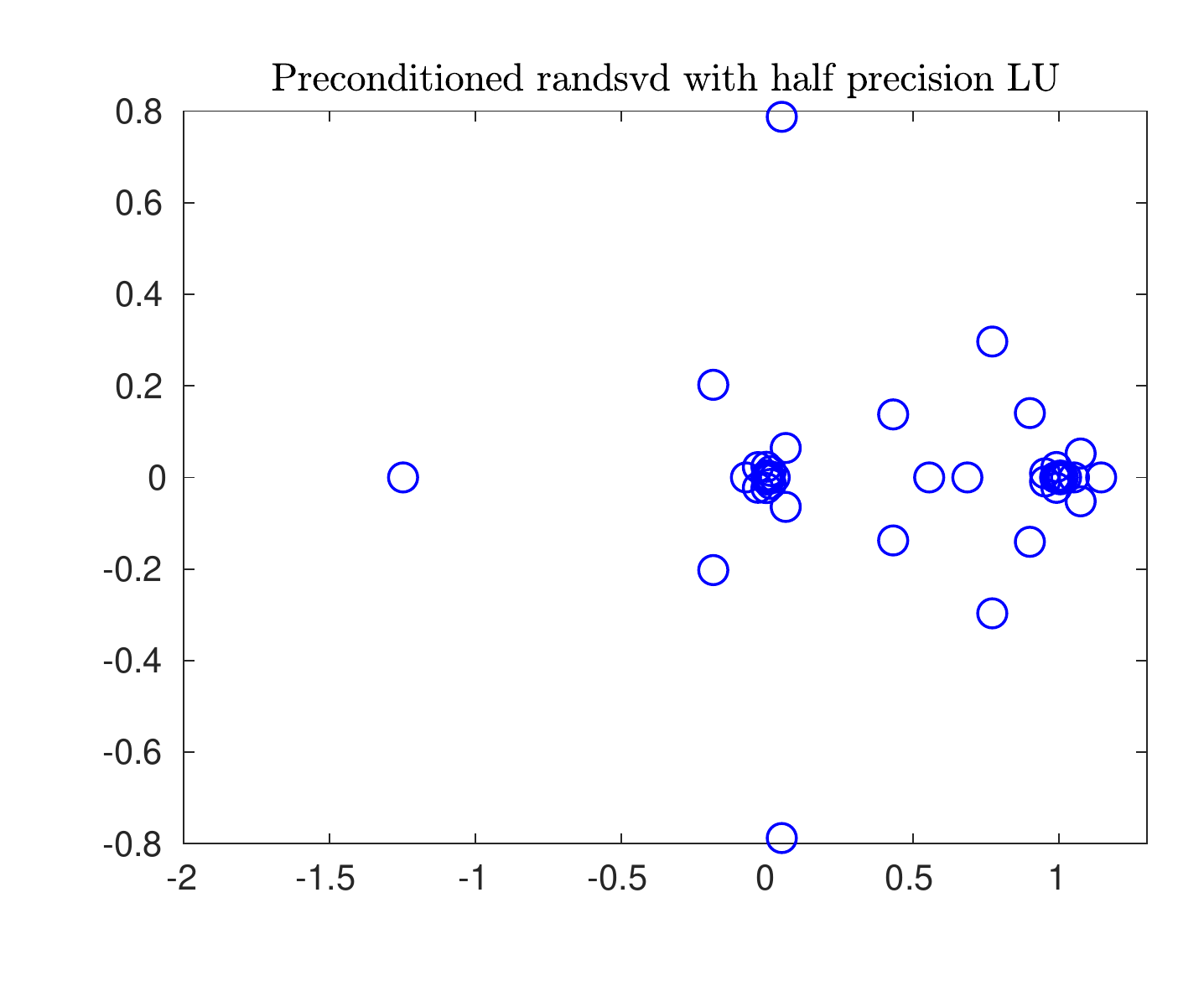}
	\caption{Eigenvalue distribution of a double-precision random dense matrix with $\kappa_2(A)=10^{12}$ without preconditioner (top left), with a double-precision LU preconditioner (top right), and with a half-precision LU preconditioner (bottom).}
	\label{fig:eig_randsvd}
\end{figure}

Figure \ref{fig:eig_randsvd} shows the eigenvalue distribution of a double-precision $100 \times 100$ random dense matrix having geometrically distributed singular values with condition number $\kappa_2(A)=10^{12}$, generated in MATLAB via the command \texttt{gallery('randsvd',100,1e12,3)}. In the unpreconditioned case (upper left), it is seen that the eigenvalues are clustered around the origin, which is a known difficult case for GMRES. When double-precision LU factors are used for preconditioning (upper right), it is observed that the eigenvalues of the preconditioned system are now clustered around 1. On the other hand, using half-precision LU factors as preconditioners (lower plot) causes a cluster of eigenvalues to remain near the origin, indicating that GMRES convergence will likely be slow (we note that these are nonnormal matrices and so the theory of \cite{lt:04} does not apply, but our experimental evidence indicates that this is the case). 


It is thus clear even in the case that low-precision LU factors can theoretically be used within GMRES-IR, they may not be the best choice from a performance perspective. In this scenario, we are left with two options: either increase the precision in which the LU factors are computed, or seek to improve the convergence behavior of GMRES through other means. It is the latter approach that we take in this work, in particular, we investigate the use of Krylov subspace recycling. 

In Section \ref{sec:recycling}, we give a background on the use of recycling in Krylov subspace methods. In Section \ref{sec:implement} we detail our implementation and experimental setup. Extensive numerical experiments that demonstrate the potential benefit of recycling within GMRES-based iterative refinement are presented in Section \ref{sec:results}. We outline future work and open problems in Section \ref{sec:conclusion}. 

\section{Krylov subspace recycling}
\label{sec:recycling}

One way to speed up the convergence of the GMRES solver is using recycling \cite{p:06, s:20}. The idea of recycling is that if we have a sequence of linear systems ($Ax_1 = b_1, Ax_2 = b_2,\ldots$) to solve and they all use the same (or a similar) coefficient matrix $A$, then we can reuse the Krylov subspace information generated in solving $Ax_1 = b_1$ to speed up converge of the method in solving $Ax_2 = b_2$, etc. This is exactly the situation we have in GMRES-IR: within the iterative refinement loop, we call GMRES on the matrix $A$ many times, and only the right-hand side changes between refinement steps. Thus we can use Krylov subspace recycling within GMRES across iterative refinement steps, and theoretically the convergence of GMRES should improve as the refinement proceeds.

GMRES-DR \cite{m:02} is a truncated and restarted solver developed for solving single nonsymmetric linear systems. The method deflates small eigenvalues for the new subspace to improve the convergence of restarted GMRES. At the end of each cycle of GMRES-DR, the desired number $k$ of harmonic Ritz values are computed, and an approximate invariant subspace associated with these values is recycled. 
One of the disadvantages of GMRES-DR is that it is not an adaptive method. 

Another truncated solver used for recycling is called GCROT \cite{d:99}. The difference between GCROT and GMRES-DR is that GCROT is a truncated minimum residual method. In other words, it recycles a subspace that minimizes the loss of orthogonality with Krylov subspace from the previous system. Like GMRES-DR, this method is also not able to be adapted for recycling.

Non-adaptive methods compute the recycling space periodically and replace it by starting from the beginning. If $A$ changes continuously, this is a problem since non-adaptive methods cannot adapt to gradual changes in matrices. 
To use truncation with adaptive recycling, the authors in \cite{p:06} modified GCROT and combined it with GMRES-DR to obtain a cheaper and more effective adaptive truncated recycling method called GCRO-DR. In this method, the residual minimization and orthogonalization are performed over the recycled subspace. This allows the algorithm to obtain sufficiently good approximations by continuously updating the approximation to the invariant subspace, and thus it requires fewer iterations per system.

The GCRO-DR algorithm uses the deflated restarting idea in GMRES-DR in the same manner as GCROT. First, using the $k$ harmonic Ritz vectors $\widetilde{Y}_k$ corresponding to the $k$ smallest harmonic Ritz values, matrices $U_k, C_k \in \mathbb{C}^{m\times k}$ are constructed such that $A^{(i+1)}U_k=C_k$ and $C_k^HC_k=I_k$ hold. After finding the optimal solution over $range(U_k)$ and computing the residual, GCRO-DR constructs the Arnoldi relation by generating a Krylov subspace of dimension $m-k+1$ with $(I-C_kC_k^H)A$. After completing the Arnoldi process, the algorithm solves a minimization problem at the end of each cycle, which reduces to an $(m+1)\times m$ least-squares problem. After solving the least-squares problem and computing the residual, a generalized eigenvalue problem is solved, and harmonic Ritz vectors are recovered. Since harmonic Ritz vectors are constructed differently than in GMRES-DR, GCRO-DR is suitable for solving individual linear systems and sequences of them. The GCRO-DR algorithm is given in Algorithm \ref{alg:gcrodr}.

\begin{algorithm}[htbp!]
	\caption{GCRO-DR \cite{p:06} \label{alg:gcrodr}}
	\begin{algorithmic}[1]
		\Require{$n\times n$ matrix $A$; right-hand-side $b$; the maximum size $m$ of the subspace; the desired number $k$ of harmonic Ritz vectors; initial guess $x_0$; GMRES convergence tolerance $\tau$.}
		\Ensure{Approximate solution $x_i$ to $Ax = b$.}
		\If {$\widetilde{Y}_k$ is known} 
		\State{Compute the reduced QR factorization $A\widetilde{Y_k}=QR$. \label{qr_cycle}} 
		\State{Compute $C_k=Q$ and $U_k = \widetilde{Y}_kR^{-1}$.}
		\State{Compute $x_1 = x_0+U_kC_k^Hr_0$ and $r_1=r_0 - C_kC_k^Hr_0$.}
		\Else
		\State{Compute $v_1=r_0/\|r_0\|_2$ and $c=\|r_0\|_2e_1$.} 
		\State{Perform $m$ steps of GMRES, solving min$\|c-\underline{H}_my\|_2$ for $y$ and generating $V_{m+1}$ and $\underline{H}_m$.}
		\State{Compute $x_1 = x_0+V_my$ and $r_1 = V_{m+1}(c-\underline{H}_my)$.}
		\State{Compute the $k$ eigenvectors $\widetilde{z}_j$ of $(H_m+h^2_{m+1,m}H^{-H}_me_me_m^H)\widetilde{z}_j = \widetilde{\theta}_j\widetilde{z}_j$ associated with the smallest magnitude eigenvalues $\widetilde{\theta}_j$ and store in $P_k$.}
		\State{Compute $\widetilde{Y}_k=V_mP_k$.}
		\State{Compute the reduced QR factorization $\underline{H}_mP_k=QR$.}
		\State{Compute $C_k = V_{m+1}Q$ and $U_k = \widetilde{Y}_kR^{-1}$.}
		\EndIf
		\While{$\|r_i\|_2>\tau$}
		\State{Perform $m-k$ Arnoldi steps with $(I-C_kC_k^H)A$, letting $v_1= r_{i-1}/\|r_{i-1}\|^2$ and generating $V_{m-k+1},\underline{H}_{m-k}$, and $B_{m-k}$.\label{arnoldi}}
		\State{Compute $\widetilde{U}_k = U_kD_k$.}
		\State{Define $\widehat{V}_m=[\widetilde{U}_k \quad V_{m-k}]$, $\widehat{W}_{m+1}=[C_k \quad V_{m-k+1}]$, and $\underline{G}_m= \left[\begin{matrix}
			D_k & B_{m-k}\\
			0 & \underline{H}_{m-k}
			\end{matrix}\right]$.}
		\State{Solve min$\|\widehat{W}_{m+1}^Hr_{i-1}-\underline{G}_my\|_2$ for $y$.}
		\State{Compute $x_i=x_{i-1}+\widehat{V}_my$ and $r_i=r_{i-1}-\widehat{W}_{m+1}\underline{G}_my$.}
		\State{Compute the $k$ eigenvectors $\widetilde{z}_i$ of $\underline{G}_m^H\underline{G}_m\widetilde{z}_i=\widetilde{\theta}_i\underline{G}_m^H\widehat{W}_{m+1}^H\widehat{V}_m\widetilde{z}_i$ associated with the smallest magnitude eigenvalues $\widetilde{\theta}_i$ and store in $P_k$.}
		\State{Compute $\widetilde{Y}_k=v_mP_k$.}
		\State{Compute the reduced QR factorization $\underline{G}_mP_k=QR$.}
		\State{Compute $C_k = W_{m+1}Q$ and $U_k = \widetilde{Y}_kR^{-1}$.}
		\EndWhile
		\State{Update $\widetilde{Y}_k = U_k$ for the next system.}
	\end{algorithmic}
\end{algorithm}

The use of recycling may also be favorable from a performance perspective. In line \ref{arnoldi} of Algorithm \ref{alg:gcrodr}, GCRO-DR performs only $m-k$ Arnoldi steps implying that it performs $m-k$ matrix-vector multiplications per cycle, whereas GMRES($m$) performs $m$ matrix-vector multiplications. It is also mentioned in \cite{p:06} that since GCRO-DR stores $U_k$ and $C_k$, it performs $2kn(1+k)$ fewer operations during the Arnoldi process. On the other hand, since we are using $k$ eigenvectors, GCRO-DR($m,k$) requires storing $k$ more vectors than GMRES($m$). Moreover, after the first cycle, in line \ref{qr_cycle}, GCRO-DR requires computing the QR factorization of $A\widetilde{Y}_k$ before the GCROT step. Although choosing the number $k$ can be considered a trade-off between memory and performance, performing $2kn(1+k)$ fewer operations in the Arnoldi stage can provide a performance improvement.


\section{Implementation and experimental setup}\label{sec:implement}

In an effort to reduce the overall computational cost of the GMRES solves within GMRES-IR, we develop a recycled GMRES-based iterative refinement algorithm, called RGMRES-IR, presented in Algorithm \ref{alg:rgmresir}. The algorithm starts with computing the LU factors of $A$ and computing the initial approximate solution in the same manner as GMRES-IR. Instead of preconditioned GMRES however, the algorithm uses preconditioned GCRO-DR($m,k$) to solve the correction equation. Similar to GMRES-IR, our presentation of RGMRES-IR in Algorithm \ref{alg:rgmresir} also uses an extra $u^2$ precision in preconditioning, although in practice one could use a uniform precision within GCRO-DR($m,k$).

\begin{algorithm}[htbp!]
	\caption{GMRES-based Iterative Refinement with Recycling (RGMRES-IR) \label{alg:rgmresir}}
	\begin{algorithmic}[1]
		\Require{$n\times n$ matrix $A$; right-hand-side $b$; maximum number of refinement steps $i_{max}$; GMRES convergence tolerance $\tau$; the maximum size $m$ of the subspace; the desired number $k$ of harmonic Ritz vectors.}
		\Ensure{Approximate solution ${x}_{i+1}$ to $Ax = b$.}
		\State{Compute LU factorization $A \approx LU$ in precision $u_f$.}
		\State{Solve $Ax_0 = b$ by substitution in precision $u_f$; store $x_0$ in precision $u$.}
		\For{$i = 0$ : $i_{max} - 1$}   
		\State{Compute $r_i = b - Ax_i$ in precision $u_r$; scale $r_i = r_i/\|r_i\|_\infty$; store in precision $u$.}
		\State{Solve $U^{-1}L^{-1}Ad_{i+1}=U^{-1}L^{-1}r_i$ by GCRO-DR($m,k$) with tolerance $\tau$ in precision $u$, with matrix-vector products with $\widetilde{A}=U^{-1}L^{-1}A$ computed in precision $u^2$; store $d_{i+1}$ in precision $u$.}
		\State{Compute $x_{i+1}=x_i+\|r_i\|_\infty d_{i+1}$ in precision $u$.}
		\IIf{converged} return $x_{i+1}$. \EndIIf
		\EndFor
	\end{algorithmic}
\end{algorithm}

In the RGMRES-IR algorithm, as in GMRES-IR, we use three precisions: $u_f$ is the factorization precision in which the factorization of $A$ is computed, $u$ is the working precision in which the input data $A$ and $b$ and the solution $x$ are stored, and $u_r$ is the precision used to compute the residuals. Again we note that we assume $u_r \leq u \leq u_f$. For the implementation, we adapted the MATLAB implementations of the GMRES-IR method developed in \cite{ch:18}, and the GCRO-DR method developed in \cite{p:06}. To simulate half-precision floating-point arithmetic, we use the \texttt{chop} library and associated functions from \cite{higham2019simulating}, available at \texttt{https://github.com/higham/chop} and \texttt{https://github.com/SrikaraPranesh/LowPrecision\_Simulation}. For single and double precision, we use the MATLAB built-in data types and to simulate quadruple precision we use the Advanpix multiprecision computing toolbox \cite{advanpix}. 

Using different precision settings results in different constraints on the condition number to guarantee convergence of the forward and backward errors. Although our experiments here use three different precisions, two precisions (only computing residuals in higher precision) or fixed (uniform) precision can also be used in the RGMRES-IR algorithm.
We also restrict ourselves to IEEE precisions (see Table \ref{tab:eps}), although we note that one could also use formats like bfloat16 \cite{bfloat16}. Table \ref{tab:prec_gmresIR} summarizes the error bounds for GMRES-IR with various choices of precisions. For convergence of GMRES-IR, $\kappa_\infty(A)$ should be less than the values shown in the fourth column. For detailed information about GMRES-IR, see \cite{ch:18}.

\begin{table}[]
	\centering
	\caption{Various IEEE precisions and their units roundoff. }
	\label{tab:eps}
	\begin{tabular}{|l|l|}
		\hline
		\multicolumn{1}{|c|}{Precision} & \multicolumn{1}{c|}{Unit Roundoff} \\ \hline
		fp16 (half) & $4.88\cdot 10^{-4}$ \\ \hline
		fp32 (single) & $5.96\cdot 10^{-8}$ \\ \hline
		fp64 (double) & $1.11\cdot 10^{-16}$ \\ \hline
		fp128 (quad) & $9.63\cdot 10^{-35}$ \\ \hline
	\end{tabular}
\end{table}

\begin{table}[h!]
	\centering
	\caption{Choices of IEEE standard precisions for three-precision GMRES-IR presented in \cite{ch:18}, and their convergence conditions.}
	\begin{tabular}{|cccc|ccc|}
		\hline
		\multirow{2}{*}{$u_f$} & \multirow{2}{*}{$u$} & \multirow{2}{*}{$u_r$}       & \multirow{2}{*}{$\kappa_\infty(A)$} & \multicolumn{2}{c}{Backward error}&                              \\
		&                    &                             &                                                  & Normwise & Componentwise & Forward error                         \\ \hline
		half                  & half             & \multicolumn{1}{c|}{single} & $10^{4}$                            & half   & half        & half  \\
		half                  & single             & \multicolumn{1}{c|}{double} & $10^{8}$                            & single   & single        & single                                \\
		half                  & double             & \multicolumn{1}{c|}{quad}   & $10^{12}$                            & double   & double        & double                                \\ 
		single                & single             & \multicolumn{1}{c|}{double} & $10^{8}$                           & single   & single        & single                                \\
		single                & double             & \multicolumn{1}{c|}{quad}   & $10^{16}$                            & double   & double        & double                                \\ \hline
	\end{tabular}
	\label{tab:prec_gmresIR}
\end{table}

We experiment with three categories of test matrices. We will first test our algorithm on random dense matrices of dimension $n=100$ having geometrically distributed singular values.  We generated the matrices using the MATLAB command \verb|gallery('randsvd',n,kappa(i),3)|, where \verb|kappa| is the array of the desired 2-norm condition numbers $\kappa_2 (A) =$\{$10^4$, $10^5$, $10^6$, $10^7$, $10^8$, $10^9$, $10^{10}$, $10^{11}$, $10^{12}$, $10^{13}$\}, and \verb|3| stands for the mode for singular value distribution of the matrices. Mode 3 corresponds to the matrix having geometrically distributed singular values. For reproducibility, we use the MATLAB command \verb|rng(1)| each time we run the algorithm. 
We will present the numerical results in Section \ref{sec:randnmat}.


As shown in Figure \ref{fig:eig_randsvd}, these matrices have eigenvalues clustered around the origin, which can be a difficult case for GMRES convergence. This class of artificial problems thus represents a good use case for the RGMRES-IR algorithm. We note however that the recycling approach can improve the convergence of GMRES even when the eigenvalues are not clustered around the origin. To illustrate this, we also test our algorithm on real matrices. All matrices shown in Table \ref{tab:suitesparse} are taken from the SuiteSparse Collection \cite{dh:11}. 
The numerical results are shown in Section \ref{sec:realmat}.

\begin{table}[]
	\centering
	\caption{Matrices from \cite{dh:11} used for numerical experiments and their properties.}
	\label{tab:suitesparse}
	\begin{tabular}{|c|c|c|c|c|c|}
		\hline
		Name       & Size & nnz  & $\kappa_\infty(A)$ & Group    & Kind                                 \\ \hline
		orsirr\_1  & 1030 & 1030 & 9.96E+04           & HB       & Computational Fluid Dynamics Problem \\
		comsol     & 1500 & 1500 & 3.42E+06           & Langemyr & Structural Problem                   \\
		circuit204 & 1020 & 1020 & 9.03E+09           & Yzhou    & Circuit Simulation Problem           \\ \hline
	\end{tabular}
\end{table}

Besides the matrices mentioned above, we also test our algorithm on prolate (symmetric, ill-conditioned Toeplitz matrices whose eigenvalues are distinct, lie in the interval $(0,1)$, and tend to cluster around 0 and 1) matrices \cite{v:93} of dimension $n=100$. We generated the matrices using the MATLAB command \verb|gallery('prolate',n,alpha)|, where \verb|alpha| is the array of the desired parameters $\alpha =$\{0.475, 0.47, 0.467, 0.455, 0.45, 0.4468, 0.44, 0.434\}. When $\alpha< 0.5$ is chosen to be small, it becomes difficult for GMRES-IR to solve the system since the eigenvalues skewed more towards zero. The properties of prolate matrices used in this study are shown in Table \ref{tab:prolate}. The numerical results are shown in Section \ref{sec:prolate}. 

\begin{table}[]
	\centering
	\caption{Prolate matrices used for numerical experiments and their properties.}
	\label{tab:prolate}
	\begin{tabular}{|c|c|c|}
		\hline
		$\alpha$      & $\kappa_\infty(A)$ & $\kappa_2(A)$ \\ \hline
		0.475  & 1.21E+06           & 3.60E+05      \\
		0.47   & 2.63E+07           & 7.60E+06      \\
		0.467  & 1.68E+08           & 4.79E+07      \\
		0.455  & 2.91E+11           & 8.04E+10      \\
		0.45   & 6.64E+12           & 1.82E+12      \\
		0.4468 & 4.98E+13           & 1.35E+13      \\
		0.44   & 3.41E+15           & 9.07E+14      \\
		0.434  & 3.44E+16           & 8.84E+15      \\ \hline
	\end{tabular}
\end{table}

For a fair comparison between GMRES-IR and RGMRES-IR, GMRES-IR is used with restart value $m$, which is the maximum size of the recycled space used in RGMRES-IR. Since the first refinement step of RGMRES-IR does not have a recycled subspace, it is the same as the first step of GMRES-IR. We thus expect a decrease in the number of GMRES iterations per refinement step starting from the second refinement step. 

The experiments are performed on a computer with Intel Core i7-9750H having 12 CPUs and 16 GB RAM with OS system Ubuntu 20.04.1. Our RGMRES-IR algorithm and
associated functions are available through the repository \verb|https://github.com/edoktay/rgmresir|, which includes
scripts for generating the data and plots in this work.

\section{Numerical experiments}\label{sec:results}
In this section, we present the numerical results comparing GMRES-IR and RGMRES-IR for solving $Ax=b$. To ensure that we fully exhibit the behavior of the methods, we set $i_{max}=10000$, which is large enough to allow all approaches that eventually converge sufficient time to do so. The GMRES convergence tolerance $\tau$, which appears both in Algorithms \ref{alg:gmresir} and \ref{alg:rgmresir}, dictates the stopping criterion for the inner GMRES iterations. The algorithm is considered converged if the relative (preconditioned) residual norm drops below $\tau$. In tests here with single working precision, we use $\tau=10^{-4}$. For double working precision, we use $\tau=10^{-8}$. The results are compared in two different metrics: The number of GMRES iterations per refinement step and the total number of GMRES iterations. The number of steps and iterations are shown in tables.

\begin{figure}[h!]
	\centering
	\includegraphics[trim={3cm 8cm 4cm 8cm},clip, width=.45\textwidth]{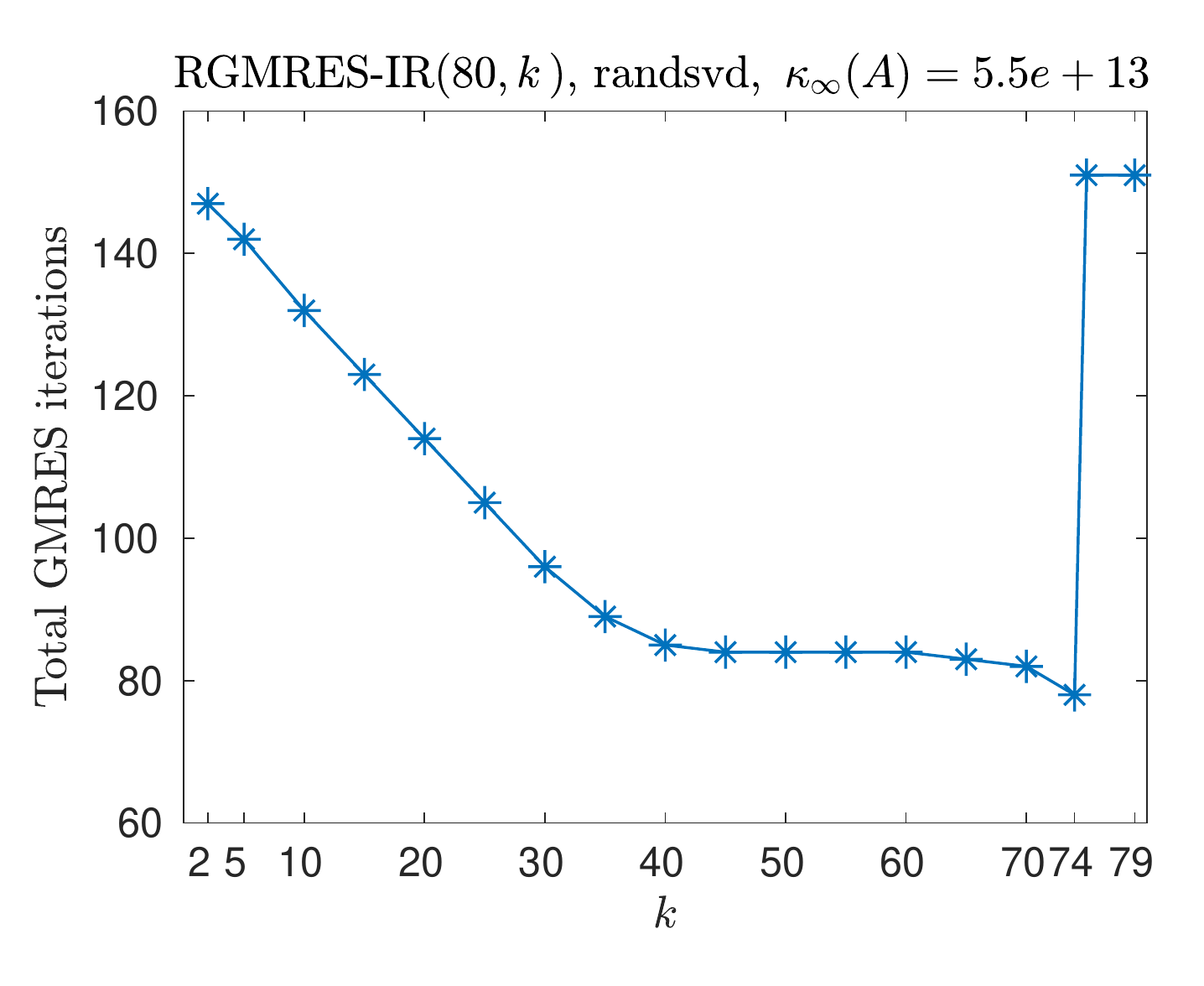}
	\includegraphics[trim={3cm 8cm 4cm 8cm},clip, width=.45\textwidth]{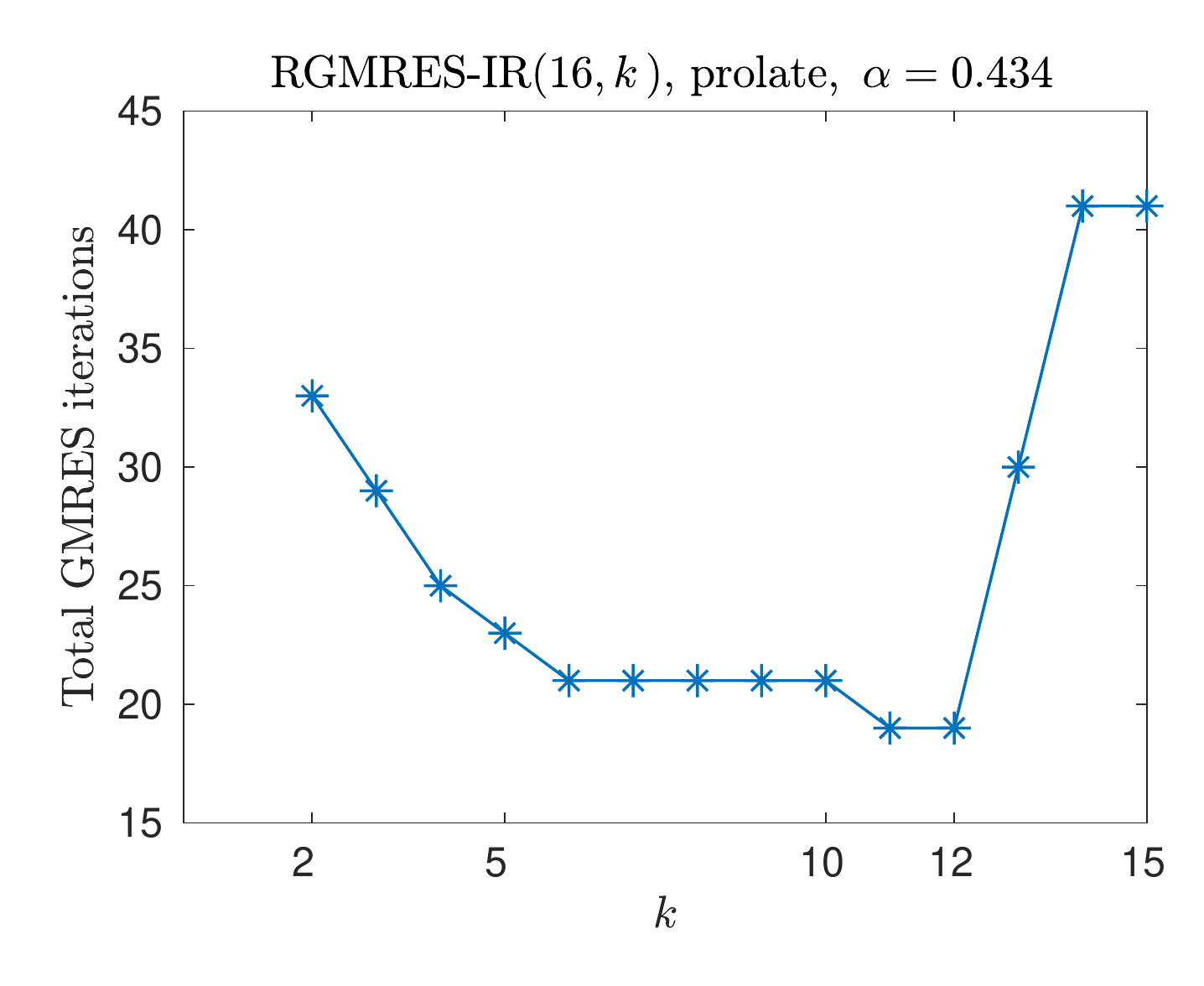}
	\caption{Total GMRES iterations for various $k$ for a randsvd matrix with $\kappa_2(A)=10^{13}$ (left) and a prolate matrix with $\alpha=0.434$ (right) for $(u_f,u,u_r)$ = (single, double, quad).}
	\label{fig:kversusgmres}
\end{figure}

In this study, $m\leq n$ is chosen such that the number of GMRES iterations in the first step of the iterative refinement method is smaller than $m$. To reduce the total number of GMRES steps in RGMRES-IR even further than the results presented in this study, one can also test RGMRES-IR with various $m\leq n$ values. However, one should note that there may be a trade-off between the total number of GMRES iterations and the number of refinement steps. For RGMRES-IR, the optimal number $k<m$ of harmonic Ritz vectors is chosen for each group of matrices with the desired precision settings after several experiments on various $(m,k)$ scenarios. The optimum $k$ differs for each matrix. The least total number of GMRES iterations is obtained for $k=(\text{the number of GMRES iterations in the first refinement step})-1$ since, in this case, we are recycling the whole generated subspace, which is expensive. That is why one should choose a $k$ value as small as possible to reduce computational cost while benefiting from recycling. Figure \ref{fig:kversusgmres} shows the change in the total GMRES iterations according to the given $k$ values for two matrices. From the plots, one can easily \textit{find the knee}, i.e., find the optimum $k$ value that uses the smallest number of GMRES iterations. 

For simplicity, $b$ is chosen to be the vector of ones for all matrices, and the precisions are chosen such that $u \leq u_f^2$, and $u_r \leq u^2$. For each random dense matrix and each matrix in Table \ref{tab:suitesparse} and Table \ref{tab:prolate}, $Ax=b$ is solved by using precisions $(u_f,u,u_r)$ = (single, double, quad), $(u_f,u,u_r)$ = (half, single, double), and $(u_f,u,u_r)$ = (half, double, quad).

In the tables we will present, the first number shows the total number of GMRES iterations. The numbers in the parentheses indicate the number of GMRES iterations performed in each refinement step. For instance, 5 (2,3) implies that at the first step of GMRES-IR, 2 GMRES iterations are performed, respectively, while for the second step, 3 GMRES iterations are performed, which gives a total of 5 GMRES iterations. 

Since all of the numerical experiments in this study were performed in MATLAB, it does not include performance analysis. Performance results for mixed precision GMRES-based iterative refinement with restarting were recently presented in \cite{lld:22}.

\subsection{Prolate matrices} \label{sec:prolate}
For the prolate matrices generated as described in Section \ref{sec:implement} with \verb|alpha| parameters $\alpha =$\{0.475, 0.47, 0.467, 0.455, 0.45, 0.4468, 0.44, 0.434\}, RGMRES-IR is used with precisions $(u_f,u,u_r)$ = (single, double, quad) and $(u_f,u,u_r)$ = (half, single, double). Tables \ref{tab:prolate_sdq} and \ref{tab:prolate_hsd} show the number of GMRES iterations performed by GMRES-IR and RGMRES-IR with different precision settings.

Table \ref{tab:prolate_sdq} shows the data for experiments with $(u_f,u,u_r)$ = (single, double, quad). In the table, we can see that for $\alpha > 0.46$, recycling does not affect the performance of GMRES for $(m,k)=(16,4)$ since only a small number of GMRES iterations are performed in each step. As $\alpha$ decreases, however, the convergence of GMRES slows down; hence recycling starts to more significantly decrease the total number of GMRES iterations. 

\begin{table}[h!]
	\centering
	\caption{Number of GMRES-IR and RGMRES-IR refinement steps and the number of GMRES iterations for each refinement step for prolate matrices with various $\alpha$ values, using precisions ($u_f,u,u_r$) = (single, double, quad) and ($m,k$) = (16,4).
	}
	\label{tab:prolate_sdq}
	\begin{tabular}{|c|cc|}
		\hline
		$\alpha$  & GMRES-IR (16)    & RGMRES-IR (16,4) \\ \hline
		0.475       & 5 (2,3) & 5 (2,3)         \\
		0.47        & 5 (2,3) & 5 (2,3)         \\
		0.467       & 7 (3,4) & 7 (3,4)         \\
		0.455       & 13 (6,7) & 8 (6,2) \\
		0.45        & 15 (7,8) & 11 (7,4)  \\
		0.4468      & 25 (7,9,9) &15 (7,4,4) \\
		0.44       & 34 (10,12,12) & 19 (10,5,4) \\
		0.434      & 41 (13,14,14) & 25 (13,6,6)  \\ \hline
	\end{tabular}
\end{table}

In the case of $(u_f,u,u_r)$ = (half, single, double), presented in Table \ref{tab:prolate_hsd}, we can see that GMRES-IR diverges for $\alpha< 0.45$ with and without recycling. However, when $\alpha=0.45$, we see that RGMRES-IR diverges although GMRES-IR converges. This is because of the multiple periods of stagnation in the second refinement step due to recycling. GMRES cannot converge in the first $16-5=11$ iterations in the second step, causing an infinite restart which results in divergence. However, for cases where both GMRES-IR and RGMRES-IR convergence, RGMRES-IR always requires fewer total GMRES iterations. For $\alpha = 0.455$, GMRES restarts in the second refinement step of GMRES-IR, while, because of recycling, RGMRES-IR converges without restarting, which decreases the computational cost.
\begin{table}[h!]
	\centering
	\caption{Number of GMRES-IR and RGMRES-IR refinement steps with the number of GMRES iterations for each refinement step for prolate matrices with various $\alpha$ values, using precisions ($u_f,u,u_r$) = (half,single,double) and ($m,k$) = (16,5).}
	\label{tab:prolate_hsd}
	\begin{tabular}{|c|cc|}
		\hline
		$\alpha$  & GMRES-IR (16) & RGMRES-IR (16,5)                                             \\ \hline
		0.475       & 12 (6,6)      & 8 (6,2)   \\
		0.47        & 16 (8,8)      & 10 (8,2) \\
		0.467       & 19 (9,10)     & 11 (9,2) \\
		0.455       & 50 (15,25,10) & 19 (15,4) \\
		0.45        & 89 (14,43,32) & -        \\
		0.4468      & -             & -          \\
		0.44        & -             & - \\
		0.434       & -             & -      \\ \hline
	\end{tabular}
\end{table}

\subsection{SuiteSparse matrices} \label{sec:realmat}
For matrices in Table \ref{tab:suitesparse}, Tables \ref{tab:real_hsd} and \ref{tab:real_hdq} compare the performance of GMRES-IR and RGMRES-IR for precisions $(u_f,u,u_r)$ = (half, single, double) and $(u_f,u,u_r)$ = (half, double, quad), respectively. In the $(u_f,u,u_r)$ = (half, single, double) scenario, for the orsirr\_1 matrix, GMRES already converges quite quickly, and so recycling does not have an effect. For the matrices comsol and circuit204, however, RGMRES-IR reduces the total number of GMRES iterations required by 35\% and 28\%, respectively. RGMRES-IR has a benefit for all matrices when using $(u_f,u,u_r)$ = (half, double, quad). 

\begin{table}[h!]
	\centering
	\caption{Number of GMRES-IR and RGMRES-IR refinement steps with the number of GMRES iterations for each refinement step for real matrices, using precisions ($u_f,u,u_r$) = (half,single,double) and ($m,k$) = (40,6).}
	\label{tab:real_hsd}
	\begin{tabular}{|c|cc|}
		\hline
		Matrix      & GMRES-IR (40) & RGMRES-IR(40,6) \\ \hline
		orsirr\_1   & 12 (6,6)      & 12 (6,6)        \\
		comsol      & 46 (22,24)    & 30 (22,8)        \\
		circuit204  & 40 (12,14,14) & 29 (12,9,8)     \\ \hline
	\end{tabular}
\end{table}

\begin{table}[h!]
	\centering
	\caption{Number of GMRES-IR and RGMRES-IR refinement steps with the number of GMRES iterations for each refinement step for real matrices, using precisions ($u_f,u,u_r$) = (half, double, quad) and ($m,k$) = (40,10).}
	\label{tab:real_hdq}
	\begin{tabular}{|c|cc|}
		\hline
		Matrix      & GMRES-IR (40)            & RGMRES-IR(40,10) \\ \hline
		orsirr\_1   & 22 (11,11) & 20 (11,9) \\
		comsol      & 52 (25,27) & 34 (25,9)   \\
		circuit204  & 59 (18,20,21) & 47 (18,14,15)  \\ \hline
	\end{tabular}
\end{table}

\subsection{Random Dense Matrices} \label{sec:randnmat}
For the random dense matrices with geometrically distributed singular values described in Section \ref{sec:implement}, we compared methods using precisions ($u_f,u,u_r$) = (single, double, quad), ($u_f,u,u_r$) = (half, single, double), and ($u_f,u,u_r$) = (half, double, quad).

\subsubsection{SGMRES-IR versus RSGMRES-IR} \label{sec:randnmat_sgmresir}

In practice, it is common that implementations use a uniform precision within GMRES (i.e., applying the preconditioned matrix to a vector in precision $u$ rather than $u^2$. This is beneficial from a performance perspective (in particular if precision $u^2$ must be implemented in software). The cost is that the constraint on condition numbers for which the refinement scheme is guaranteed to converge becomes tighter. To illustrate the benefit of recycling in this scenario, we first compare what we call SGMRES-IR (which is GMRES-IR but with a uniform precision within GMRES) to the recycled version, RSGMRES-IR. For a fair comparison, restarted SGMRES-IR (SGMRES-IR($m$)) is compared with recycled SGMRES-IR (RSGMRES-IR($m,k$)).

Table \ref{tab:sgmresir} shows the number of GMRES iterations performed by SGMRES-IR and RSGMRES-IR in the ($u_f,u,u_r$) = (single, double, quad) setting.
From the table, we observe that recycling reduces the number of GMRES iterations in this case as well. The reason why SGMRES-IR does not converge for $\kappa_2(A)\geq10^{14}$ is that in the first refinement step, restarted SGMRES does not converge (restarting an infinite number of times). For RSGMRES-IR, in the first GCRO-DR call, the recycling after the first restart cycle helps, allowing GCRO-DR to converge. We note that this is another benefit of the recycling approach, as it can improve the reliability of restarted GMRES, which is almost always used in practice. 

\begin{table}[h!]
	\centering
	\caption{Number of SGMRES-IR and RSGMRES-IR refinement steps with the number of GMRES iterations for each refinement step for precisions ($u_f,u,u_r$) = (single, double, quad) and ($m,k$) = (80,18). For $\kappa_2(A)=10^{15}$, RSGMRES-IR required 2093 (139, 60, 59, 60, 59, 59, 60, 60, 60, 60, 60, 60, 60, 60, 60, 60, 60, 48, 46, 46, 51, 60, 48, 48, 48, 51, 48, 48, 48, 58, 48, 48, 48, 48, 48, 48, 61) iterations.}
	\label{tab:sgmresir}
	\begin{tabular}{|cc|cc|}
		\hline
		$\kappa_\infty(A)$ & $\kappa_2(A)$  & SGMRES-IR (80)      & RSGMRES-IR (80,18) \\ \hline
		$6 \cdot 10^9$     & $10^9$          & 64 (19,23,22) & 34 (19,8,7) \\
		$6 \cdot 10^{10}$  & $10^{10}$       & 120 (39,40,41) & 65 (39,13,13) \\
		$6 \cdot 10^{11}$  & $10^{11}$       & 160 (52,54,54) & 94 (52,21,21) 	\\
		$6 \cdot 10^{12}$  & $10^{12}$       & 196 (65,65,66) & 163 (65,32,32,34) \\
		$5 \cdot 10^{13}$  & $10^{13}$       & 301 (75,75,75,76) & 199 (75,41,41,42) \\
		$5 \cdot 10^{14}$  & $10^{14}$       & - & 493 (131,51,51,52,52,52,52,52) \\
		$5 \cdot 10^{15}$  & $10^{15}$       & - & 2093* \\ \hline
	\end{tabular}
\end{table}

\subsubsection{GMRES-IR versus RGMRES-IR} \label{sec:randnmat_gmresir}
We now return to our usual setting and compare GMRES-IR and RGMRES-IR for random dense matrices with condition numbers $\kappa_2 (A) =$ \{$10^4$, $10^5$, $10^6$, $10^7$, $10^8$, $10^9$, $10^{10}$, $10^{11}$, $10^{12}$, $10^{13}$, $10^{14}$, $10^{14}$\}. Results using precisions $(u_f,u,u_r)$ = (single, double, quad), $(u_f,u,u_r)$ = (half, single, double), and $(u_f,u,u_r)$ = (half, double, quad) are displayed in Tables \ref{tab:randn_sdq}-\ref{tab:randn_hdq}, respectively.

Table \ref{tab:randn_sdq} shows the numerical experiments with $(u_f,u,u_r)$ = (single, double, quad). In the table, we can see that for relatively well-conditioned matrices ($\kappa_\infty(A)<10^{8}$) since the total number of GMRES iterations is small, using recycling does not change the total number of GMRES iterations. When the condition number increases, however, the number of GMRES iterations drops significantly. For $\kappa_2(A)\geq10^{14}$, GMRES-IR does not converge for the same reason SGMRES-IR did not as described above; namely, that GMRES with the restart parameter 80 does not converge. As before, the recycling between restart cycles fixed this, and RGMRES-IR is still able to converge in this case. 

\begin{table}[h!]
	\centering
	\caption{Number of GMRES-IR and RGMRES-IR refinement steps with the number of GMRES iterations for each refinement step for random dense matrices having geometrically distributed singular values (mode 3) with various condition numbers, using precisions ($u_f,u,u_r$) = (single, double, quad)  and ($m,k$) = (80,18).}
	\label{tab:randn_sdq}
	\begin{tabular}{|cc|cc|}
		\hline
		$\kappa_\infty(A)$ & $\kappa_2(A)$  & GMRES-IR (80)       & RGMRES-IR ($80,18$) \\ \hline
		$9 \cdot 10^4$     & $10^4$         & 4 (2,2) & 4 (2,2) \\
		$8 \cdot 10^5$     & $10^5$         & 6 (3,3) & 6 (3,3) \\
		$7 \cdot 10^6$     & $10^6$         & 8 (4,4) & 8 (4,4) \\
		$7 \cdot 10^7$     & $10^7$         & 11 (5,6) & 11 (5,6) \\
		$7 \cdot 10^8$     & $10^8$         & 22 (10,12) & 22 (10,12) \\
		$6 \cdot 10^9$     & $10^9$         & 67 (19,24,24) & 36 (19,8,9) \\
		$6 \cdot 10^{10}$  & $10^{10}$      & 80 (39,41) & 53 (39,14) \\
		$6 \cdot 10^{11}$  & $10^{11}$      & 107 (52,55) & 75 (52,23) \\
		$6 \cdot 10^{12}$  & $10^{12}$      & 132 (65,67) & 99 (65,34) \\
		$5 \cdot 10^{13}$  & $10^{13}$      & 151 (75,76) & 117 (75,42) \\
		$5 \cdot 10^{14}$  & $10^{14}$      & - & 184 (131,53) \\
		$5 \cdot 10^{15}$  & $10^{15}$      & - & 325 (139,62,124) \\ \hline
	\end{tabular}
\end{table}

In Table \ref{tab:randn_hsd} we present the experiments for ($u_f,u,u_r$) = (half, single, double). For $\kappa_\infty(A)>10^5$, recycling reduces the total number of GMRES iterations. This is also the case in Table \ref{tab:randn_hdq}, which shows results for precisions ($u_f,u,u_r$) = (half, double, quad). This is a known difficult case for GMRES, and thus is a clear case where we can see significant benefit of recycling. We see the most significant improvement for the matrix with $\kappa_\infty(A)=10^{13}$, in which RGMRES-IR requires over $16\times$ fewer GMRES iterations than GMRES-IR. We note that GMRES-IR is only guaranteed to converge up to $\kappa_2(A)<10^{12}$ for this combination of precisions; see Table \ref{tab:prec_gmresIR}.

\begin{table}[h!]
	\centering
	\caption{Number of GMRES-IR and RGMRES-IR refinement steps with the number of GMRES iterations for each refinement step for random dense matrices having geometrically distributed singular values (mode 3) with various condition numbers, using precisions ($u_f,u,u_r$) = (half, single, double)  and ($m,k$) = (90,30).}
	\label{tab:randn_hsd}
	\begin{tabular}{|cc|cc|}
		\hline
		$\kappa_\infty(A)$ & $\kappa_2(A)$   & GMRES-IR (90) & RGMRES-IR ($90,30$) \\ \hline
		$9 \cdot 10^4$     & $10^4$   & 20 (9,11)     & 20 (9,11)           \\
		$8 \cdot 10^5$     & $10^5$  & 70 (32,38)    & 44 (32,12)         \\
		$7 \cdot 10^6$     & $10^6$  & 129 (64,65)   & 77 (64,13)
		\\
		$7 \cdot 10^7$     & $10^7$  & 164 (82,82)   & 107 (82,25)
		\\ \hline
	\end{tabular}
\end{table}

\begin{table}[h!]
	\centering
	\caption{Number of GMRES-IR and RGMRES-IR refinement steps with the number of GMRES iterations for each refinement step for random dense matrices having geometrically distributed singular values (mode 3) with various condition numbers, using precisions ($u_f,u,u_r$) = (half, double, quad)  and ($m,k$) = (100,30). For $10^{13}$, GMRES-IR did 3954 (102, 105, 98, 98, 88, 81, 85, 84, 81, 84, 83, 84, 89, 90, 87, 88, 89, 91, 83, 83, 87, 88, 91, 84, 89, 88, 83, 82, 89, 92, 88, 97, 83, 84, 90, 83, 84, 83, 90, 83, 94, 83, 82, 85, 99) iterations.}
	\label{tab:randn_hdq}
	\begin{tabular}{|cc|cc|}
		\hline
		$\kappa_\infty(A)$ & $\kappa_2(A)$  & GMRES-IR (100)       & RGMRES-IR (100,30) \\ \hline
		$9 \cdot 10^4$     & $10^4$         & 33 (16,17) & 33 (16,17) \\
		$8 \cdot 10^5$     & $10^5$         & 85 (41,44) & 71 (41,15,15) \\
		$7 \cdot 10^6$     & $10^6$        & 134 (66,68) & 85 (66,19) \\
		$7 \cdot 10^7$     & $10^7$         & 167 (83,84)  & 113 (83,30) \\
		$7 \cdot 10^8$     & $10^8$         & 193 (96,97) & 138 (96,42) \\
		$6 \cdot 10^9$     & $10^9$         & 200 (100,100)  & 151 (100,51) \\
		$6 \cdot 10^{10}$  & $10^{10}$      & 200 (100,100)        & 158 (100,58) \\
		$6 \cdot 10^{11}$  & $10^{11}$      & 200 (100,100)        & 165 (100,65) \\
		$6 \cdot 10^{12}$  & $10^{12}$      & 200 (100,100)        & 170 (100,70) \\
		$5 \cdot 10^{13}$  & $10^{13}$       & 3954*    & 241 (171,70) \\ \hline
	\end{tabular}
\end{table}

\begin{figure}[h!]
	\centering
	\includegraphics[trim={3cm 8cm 4cm 8cm},clip, width=.45\textwidth]{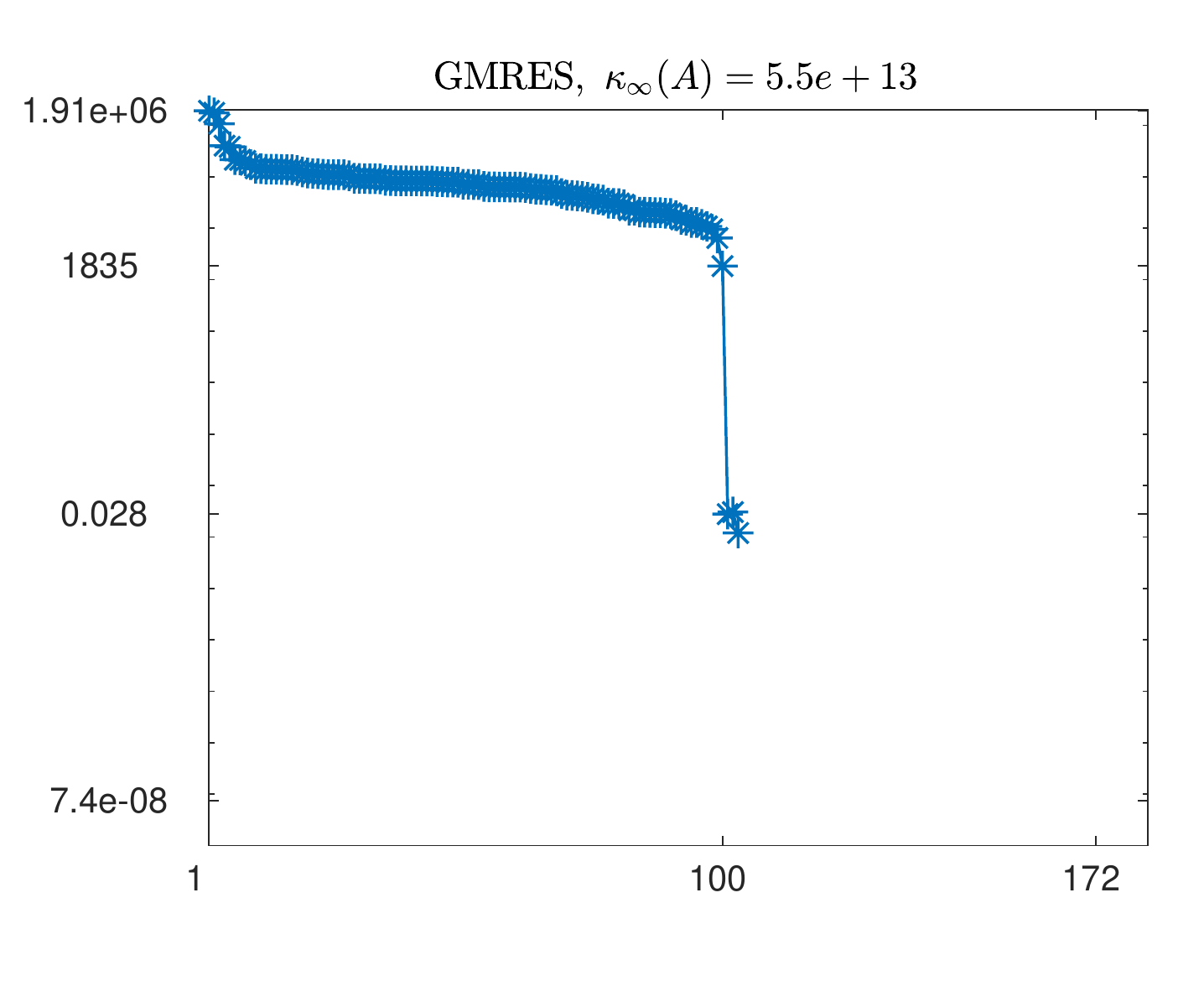}
	\includegraphics[trim={3cm 8cm 4cm 8cm},clip, width=.45\textwidth]{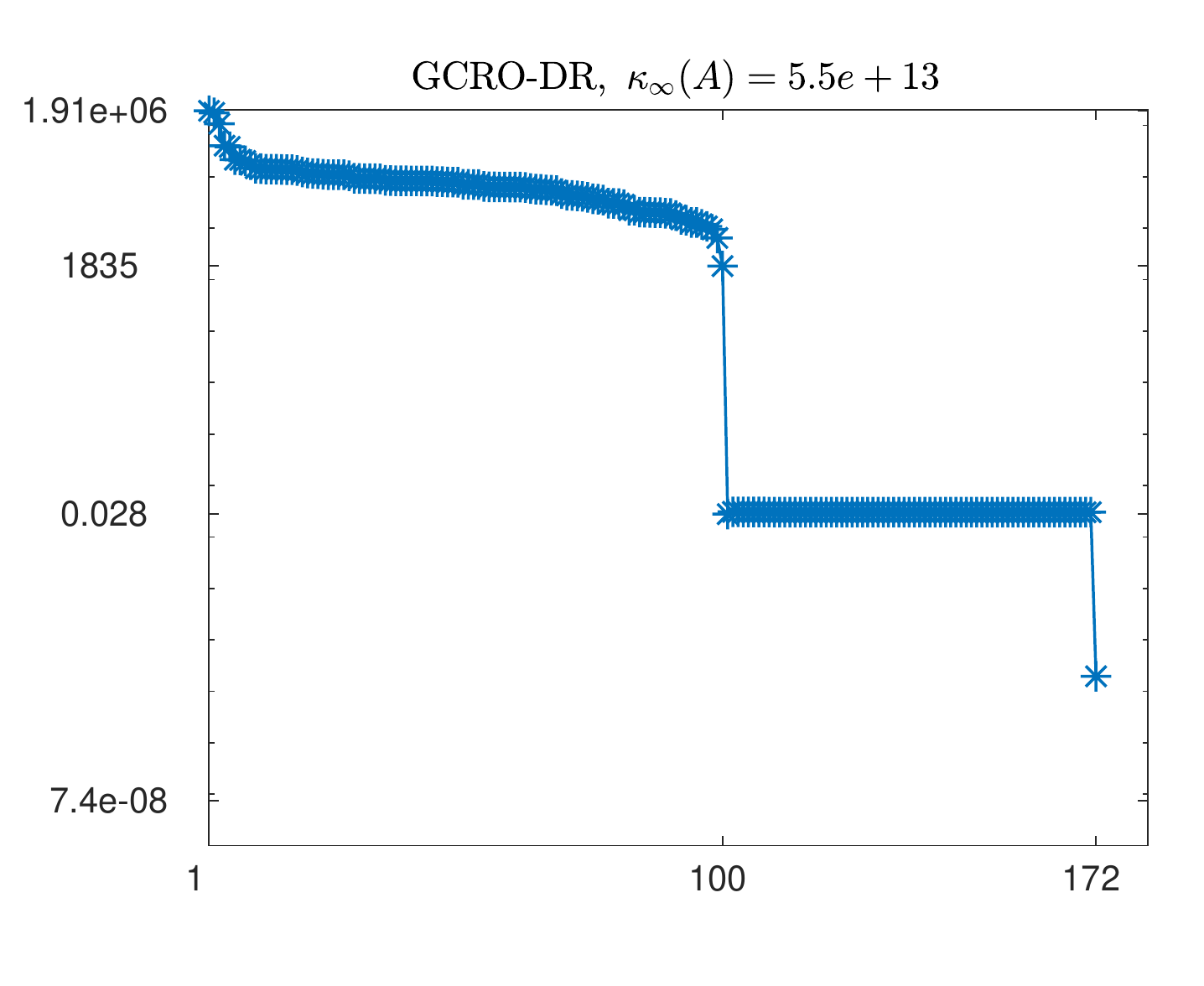}
	\caption{Residual trajectory of GMRES (left) and GCRO-DR (right), used within GMRES-IR and RGMRES-IR, respectively, for a randsvd matrix with $\kappa_2(A)=10^{13}$ and precisions $(u_f,u,u_r)$ = (half, double, quad).}
	\label{fig:error}
\end{figure}

The reason that RGMRES-IR outperforms GMRES-IR is different than in the previous cases (caused by stagnation caused by restarting), and is almost accidental in this case. We investigate this more closely in Figure \ref{fig:error}. In the left plot, we see the convergence trajectory of GMRES(100). In the first restart cycle, the residual decreases from $10^{6}$ to $10^{3}$ after 100 GMRES iterations. GMRES restarts and performs two more iterations, at which point it converges to a relative residual of $10^{-8}$ (absolute residual of around $10^{-2}$. Hence, the first refinement step of GMRES-IR does $100+2=102$ iterations. The right plot shows the residual trajectory for GCRO-DR. The first restart cycle is the same as in GMRES, however, once the method restarts, the residual stagnates just above the level required to declare convergence. After $m-k=70$ more iterations, GCRO-DR restarts again, and this time, the residual drops significantly lower. So while GCRO-DR requires more iterations (172) to converge to the specified tolerance, when it does converge, it converges to a solution with a smaller residual. This phenomenon can in turn reduce the total number of refinement steps required. It is possible that we could reduce the overall number of GMRES iterations within GMRES-IR (and also RGMRES-IR) by making the GMRES convergence tolerance $\tau$ smaller. We did not experiment with changing the GMRES tolerance within GMRES-IR or RGMRES-IR, but this tradeoff would be interesting to explore in the future. 


We stress that the convergence guarantees for GMRES-IR for various precisions stated in \cite{ch:17, ch:18, h:21} hold only for the case of unrestarted GMRES, i.e., $m=n$. When $m<n$, there is no guarantee that GMRES converges to a backward stable solution and thus no guarantee that GMRES-IR will converge. Choosing a restart parameter $m$ that allows for convergence is a difficult problem, and a full theory regarding the behavior of restarted GMRES is lacking. The behavior of restarted GMRES is often unintuitive; whereas one would think that a larger restart parameter is likely to be better than a smaller one as it is closer to unrestarted GMRES, this is not always the case. In \cite{e:03}, the author gave examples where a larger restart parameter causes complete stagnation, whereas a smaller one results in fast convergence.

Whereas in Table \ref{tab:randn_hdq_90} we now illustrates the behavior for the same problems and same precisions ($u_f,u,u_r$) = (half, double, quad) for the case $m<n$. It is seen from the table that both methods converge for $\kappa_\infty(A)<10^8$. After this point, GMRES-IR does not converge, whereas RGMRES-IR does. This serves as an example where the convergence guarantees given in \cite{ch:17,ch:18} do not hold for GMRES-IR with restarted GMRES; for unrestarted GMRES, convergence is guaranteed up to $\kappa_\infty(A)\leq 10^{12}$ for this precision setting. Here, GMRES-IR does not converge because of the stagnation caused by restarting in the first refinement step. Aided by the recycling between restart cycles, RGMRES-IR does converge up to $\kappa_2(A)=10^{11}$, although the large number of GMRES iterations required in the first refinement step makes this approach impractical. 

\begin{table}[]
	\centering
	\caption{Number of GMRES-IR and RGMRES-IR refinement steps with the number of GMRES iterations for each refinement step for random dense matrices having geometrically distributed singular values (mode 3) with various condition numbers, using precisions ($u_f,u,u_r$) = (half, double, quad)  and ($m,k$) = (90,40).}
	\label{tab:randn_hdq_90}
	\begin{tabular}{|cc|cc|}
		\hline
		$\kappa_\infty(A)$ & $\kappa_2(A)$ & GMRES-IR (90) & RGMRES-IR (90,40) \\ \hline
		$9 \cdot 10^4$     & $10^4$        & 33 (16,17)    & 33 (16,17)        \\
		$8 \cdot 10^5$     & $10^5$        & 85 (41,44)    & 50 (41,9)         \\
		$7 \cdot 10^6$     & $10^6$        & 134 (66,68)   & 81 (66,15)        \\
		$7 \cdot 10^7$     & $10^7$        & 167 (83,84)   & 100 (83,17)       \\
		$7 \cdot 10^8$     & $10^8$        & -             & 149 (119,30)      \\
		$6 \cdot 10^9$     & $10^9$        & -             & 179 (134,45)      \\
		$6 \cdot 10^{10}$  & $10^{10}$     & -             & 470 (388,41,41)   \\
		$6 \cdot 10^{11}$  & $10^{11}$     & -             & -                 \\
		$6 \cdot 10^{12}$  & $10^{12}$     & -             & -                 \\
		$5 \cdot 10^{13}$  & $10^{13}$     & -             & -                 \\ \hline
	\end{tabular}
\end{table}



\section{Conclusion and future work}
\label{sec:conclusion}
With the emergence of mixed precision hardware, mixed precision iterative refinement algorithms are the focus of significant renewed interest. A promising approach is the class of GMRES-based refinement schemes, which can enable the accurate solution of extremely ill-conditioned matrices. However, for some matrices, GMRES convergence can be very slow, even when (low-precision) preconditioners are applied. This makes the GMRES-based approaches unattractive from a performance perspective. In this work, incorporate Krylov subspace recycling into the mixed precision GMRES-based iterative refinement algorithm in order to reduce the total number of GMRES iterations required. We call our algorithm RGMRES-IR. Instead of preconditioned GMRES, RGMRES-IR uses a preconditioned GCRO-DR algorithm to solve for the approximate solution update in each refinement step.  Our detailed numerical experiments on random dense matrices, prolate matrices, and matrices from SuiteSparse \cite{dh:11} show the potential benefit of the recycling approach. Even in cases where the number of GMRES iterations does not preclude the use of GMRES-based iterative refinement, recycling can have a benefit. In particular, it can improve the reliability of restarted GMRES, which is used in most practical scenarios. 

One major caveat for GMRES-based iterative refinement schemes is that the analysis and convergence criteria discussed in the literature all rely on the use of unrestarted GMRES. When restarted GMRES is used, we can not give such concrete guarantees, as restarted GMRES may not converge even in infinite precision. A greater understanding of the theoretical behavior of restarted GMRES (and GCRO-DR) both in infinite and finite precision would be of great interest. 

Another potential future direction is the exploration of the potential for the use of mixed precision within GCRO-DR. In this work, within GCRO-DR we only used extra precision in applying the preconditioned matrix to a vector, as is done in GMRES-IR. There may however, be further potential for the use of low precision within GCRO-DR, for example, in computation of harmonic Ritz pairs. 

\clearpage
\bibliographystyle{siamplain}
\bibliography{paper}

\end{document}